\documentclass[11pt]{article}
\usepackage{amssymb}

\begin{document}

\begin{center}
\textbf{MEAN CURVATURE MOTION OF TRIPLE JUNCTIONS OF GRAPHS IN TWO DIMENSIONS}\\
by A. Freire\\

University of Tennessee, Knoxville \end{center}

\begin{abstract}
We consider a system of three surfaces, graphs over a bounded domain
in ${\mathbb R}^2$, intersecting along a time-dependent curve and
moving by mean  curvature while preserving the pairwise angles at
the curve of intersection (equal to $2\pi/3$.) For the corresponding
two-dimensional parabolic free boundary problem we prove short-time
existence of classical solutions (in parabolic H\"{o}lder spaces),
for sufficiently regular initial data satisfying a compatibility
condition.
\end{abstract}

\vspace{.5cm}

\textbf{0. Introduction: Triple junctions of graphs.}\vspace{.2cm}

The goal of this paper is to prove short-time existence for the
following geometric evolution problem: a configuration of three
surfaces in $\mathbb{R}^3$, bordered by a common time-dependent
curve of intersection $\Lambda(t)$, moves with normal velocity given
by mean curvature (at each interior point) in such a way that the
pairwise angles defined by their unit normals along $\Lambda(t)$ are
constant throughout the evolution, equal to $2\pi/3$
radians.\vspace{.2cm}

The corresponding evolution for systems of curves was considered in
\cite{BronsardReitich} (where short-time existence was proved) and
in \cite{Mantegazza et al.}, which includes a continuation criterion
and the blow-up analysis at the first singular time. Very recently,
the preprint \cite{LensSeminar} includes global existence results
for symmetric systems of curves of `lens type'. In addition to being
the natural time-dependent version of the classical problem of
minimal surfaces meeting along a `liquid edge' (\cite{Dierkes et
al.} p. 299 ff.), the problem derives its interest from being the
sharp-interface limit of a well-known class of parabolic
evolutionary models in materials science, defined by a vector-valued
order parameter and a multi-well potential. (See
\cite{BronsardReitich} for a discussion in the case of
curves.)\vspace{.2cm}

In the case of surfaces, the parametric approach adopted to prove
local existence for curves in \cite{BronsardReitich} and
\cite{Mantegazza et al.} does not work, that is, does not lead
directly to a well-posed parabolic system. The essential difficulty
can be traced back to the fact that the junction is now a
one-dimensional object, and one does not expect that independently
evolving parametrizations of each surface in a triple junction
configuration will continue to match pointwise along the junction,
even if they do so at $t=0$. On the other hand, it is difficult to
encode the matching condition analytically so as to allow for such
`sliding' along the junction. Thus one is naturally led to consider
the special case of \emph{graphs}, for the simple reason that the
graph parametrization is canonical, so that the surfaces' meeting
along $\Lambda(t)$ corresponds to three functions coinciding
pointwise along a moving interface $\Gamma_t$ in $\mathbb{R}^2$.
\vspace{.2cm}

In this way the problem for graphs becomes a more-or-less classical
free boundary problem for a quasilinear parabolic system in two
dimensions. By analogy with the usual Stefan problem terminology,
this would be a `three-phase problem', with two of the `phases'
coexisting in the same domain. The corresponding `one-phase'
problem- the mean curvature motion of a single graph over a
time-dependent domain $D(t)$ in $\mathbb{R}^n$, intersecting
$\mathbb{R}^n$ at a constant angle- has been dealt with in a recent
preprint of the author (\cite{Freire}), including local existence
and some results on long-time behavior. In \cite{Freire}, we
introduced additional `orthogonality' conditions at the free
boundary to obtain a well-posed system. If one tries to extend this
method to the more general case considered here, compositions of the
unknown functions with the (also unknown) diffeomorphisms pulling
back the evolution to a fixed domain come into play, and since
compositions behave poorly in H\"{o}lder spaces it becomes difficult
to implement a fixed-point scheme using this method. Instead, we use
here an ingenious transformation method which we learned from the
paper \cite{BaconneauLunardi}. (These authors considered a one-phase
multidimensional free-boundary problem for a semilinear parabolic
system).\vspace{.2cm}

Now for a precise description of the problem and results. We focus
in this paper on configurations parametrized by two disks and one
annulus, but it is clear the method would work for more general
configurations (for example, two annuli over one sub-domain, one
disk over the other; or three annuli.) In addition, the system lives
in a fixed bounded domain $\Omega \subset {\mathbb R^2}$, and we
must specify how one of the surfaces intersects $\partial\Omega$; we
choose orthogonal intersection, which leads to a standard Neumann
condition at the fixed boundary. More generally, one could prescribe
orthogonal intersection with a non-vertical, fixed `support
surface'. This would introduce additional moving boundaries into the
problem but analytically would present no essential difficulty, so
we avoided this in the interest of simplicity.\vspace{.2cm}

A \emph{configuration} at time $t$  consists of three surfaces:
$$\Sigma^I_t=graph(w^I_t),\quad w^I_t=w^I(\cdot,t):D^I_t\rightarrow
\mathbb{R},\quad I=1,2,3,$$ where $D^1_t=D^2_t$ and
$D^3_t=\Omega-\bar{D}^1_t$ are bounded, connected open subsets of
the bounded domain $\Omega\subset \mathbb{R}^2$. The $\Sigma^I_t$
intersect along an embedded closed curve $\Lambda_t\subset
\mathbb{R}^3$, a graph over the simple closed planar curve
$\Gamma_t=\bar{D}^1_t\cap\bar{D}^3_t$. The upper unit normals are:
$$N^I=\frac 1{v^I}[-Dw^I,1],\quad v^I=\sqrt{1+|Dw^I|^2}.$$
Along $\Gamma_t$, we impose the following conditions:

(i) (matching) $w^1(y,t)=w^2(y,t)=w^3(y,t),\quad y\in \Gamma_t\quad
(BC^0_w)$.

(ii) (angle) $N^1+N^2=N^3$.

\noindent That is, the $\Sigma^I$ intersect making pairwise $2\pi/3$
angles. In addition, $\Sigma^3$ intersects the vertical cylinder
over $\partial \Omega$ orthogonally: if $\omega$ denotes the inner
unit normal to $\partial \Omega$, we require:
$$d_{\omega}w^3=0\mbox{ on }\partial \Omega.\quad (BC_{w^3})$$

The angle condition is equivalent to two scalar conditions on
$\Gamma_t$:
$$\left \{ \begin{array}{lr}
\frac 1{v^1}+\frac 1{v^2}=\frac 1{v^3}&(BC_w^1)\\
\frac{d_nw^1}{v^1}+\frac{d_nw^2}{v^2}=\frac{d_nw^3}{v^3},&(BC_w^2)
\end{array}\right .$$
where $n$ is the inner unit normal of $\Gamma_t$.

\vspace{.2cm}

To specify a configuration, we need three time-dependent functions
$w^I(t):D^I_t\rightarrow \mathbb R$, defining graph parametrizations
of the $\Sigma^I_t$:
$$G^I:E^I_T\rightarrow \mathbb R^3, \quad E_T^I=\{(x,t)\in \mathbb{R}^2\times
[0,T] ;x\in D^I_t\}$$
$$G^I(x,t)=[x,w^I(x,t)]\in {\mathbb R}^2\times {\mathbb R}.$$
(The labeling of the surfaces is chosen so that near $\Lambda(t)$,
$\Sigma^2$ lies below $\Sigma^1$.)\vspace{.2cm}

We wish to evolve a given initial configuration through
configurations, so that on each surface the normal velocity is given
by the mean curvature at each point:
$$\langle \partial_tG^I,N^I\rangle =H^I.$$
Recall $H^I=tr_{g^I}h^I=g^{Iij}h^I_{ij}$, where $h^I$ is the second
fundamental form of $\Sigma^I$, pulled back to $\mathbb{R}^2$ via
$G$:
$$h^I(e,e)=-\langle d_eN^I,d_eG^I\rangle=\langle
d^2G^I(e,e),N^I\rangle=\frac 1{v^I}d^2w^I(e,e),\quad e\in
\mathbb{R}^2,$$ and $g^I$ is the induced metric:
$$g^I_{ij}=\langle
\partial_iG^I,\partial_jG^I\rangle=\delta_{ij}+\partial_iw^I\partial_jw^I.$$
Mean curvature motion of graphs is given by the equation:
$$(MCM)_w\qquad L_{g^I}[w^I]=0,\quad
L_{g^I}[f]:=\partial_tf-tr_{g^I}D^2f.$$\vspace{.2cm}

Our local existence result follows:\vspace{.2cm}

\textbf{Theorem.} \emph{Let $\Omega\subset {\mathbb R}^2$ be a
bounded open set with $C^{3+\alpha}$ boundary, $D_0^1=D_0^2$ and
$D_0^3=\Omega-\bar{D}_0^1$ open subsets of $\Omega$ with common
boundary component $\Gamma_0$, a simple closed $C^{3+\alpha}$ curve
in $\Omega$. Given $w_0^I\in C^{3+\alpha}(D_0^I)$, assume}:

\emph{(1) The $\Sigma_0^I=graph(w_0^I)$ define a triple junction
configuration; }

 \emph{(2) The mean curvatures $H_0^I$ of $\Sigma_0^I$
satisfy at points of $\Gamma_0$ the following compatibility
condition:}
$$H_0^1+H_0^2=H_0^3,\quad y\in \Gamma_0.$$

 \emph{Then there exists $T>0$ depending only on the initial data, domains
$E_T^I\subset \Omega\times [0,T]$ and functions $w^I\in
C^{2+\alpha,1+\alpha/2}(E_T^I)$ such that:}

\emph{(3) $w^I$ is a classical solution of $L_{g^I}[w^I]=0$ in
$E^I$; }

 \emph{(4) The conjugation conditions $(BC_w)$ hold on
$\Gamma_t=\partial D^1_t$, and the boundary condition $(BC_{w^3})$
holds on $\partial \Omega$.}\vspace{.4cm}

\emph{Notation.} $$E_T^I=\bigcup_{t\in [0,T]}D^I_t\times\{t\}\subset
\Omega \times [0,T];$$ it is always assumed that $T=\sup\{t\in
[0,T); D^I_t\neq \emptyset \forall I\}$. We use standard notation
for parabolic H\"{o}lder spaces; for example, the norm in
$C^{2+\alpha,1+\alpha/2}(E)$ (where $E\subset \Omega \times [0,T]$)
is given by:
$$||f||_{2+\alpha}:=||f||_{C^{2,1}(E)}+[D^2f]_x^{(\alpha)}+[D^2f]_t^{(\alpha/2)}+[\partial_tf]_x^{(\alpha)}+[\partial_t
f]_t^{(\alpha/2)}+[Df]_t^{((1+\alpha)/2)},$$ with the standard
notations for H\"{o}lder-type difference quotients in $x$ or in $t$.
$\alpha$ is an arbitrary constant in the open interval $(0,1)$,
fixed throughout the paper.

The subscript $t$ always denotes a map (or function, or set) at time
$t$, never a partial derivative (denoted $\partial_t$). The unit
vectors $\nu$, $\tau_0$ are the inner unit normal (resp.
counterclockwise unit tangent) vectors to $\Gamma_0$, while $n$ and
$\tau$ denote the corresponding unit vectors for
$\Gamma_t$.\vspace{.3cm}

\emph{Remarks.}

(1) Note the drop in regularity: we require $C^{3+\alpha}$ initial
data for $C^{2+\alpha}$ (in space) solutions. This occurs also in
the main result in \cite{BaconneauLunardi}, and for the same reason.
It was also observed in \cite{Freire}, where local existence for the
`one-phase' problem is obtained by a different
argument.\vspace{.2cm}

(2) The result, at this point, lacks both a uniqueness statement and
a continuation criterion. It is conceivable that assuming higher
regularity on the initial data ($C^{4+\alpha}$), a uniqueness
statement could be obtained following the argument in
\cite{BaconneauLunardi}. Both issues will be considered in further
work.\vspace{.2cm}

(3) The restriction to two dimensions makes the argument technically
easier, but is probably inessential. On the other hand, generalizing
to configurations beyond graphs may require a completely different
line of argument.\vspace{.3cm}

\emph{Outline.} Section 1: We derive the compatibility condition and
the equation of motion for the junction $\Lambda(t)$ (and for the
`interface' $\Gamma_t$). Section 2:  In a standard way, we use a
time-dependent diffeomorphism to pass to a system (with \emph{four}
unknown functions) over the fixed domains $D_0^I$:
$$w^I(y,t),\quad y\in D^I_t\longrightarrow u^I(x,t),\tilde{\rho}(x,t),\quad x\in D_0^I,$$
where $\tilde{\rho}$ is an extension to $\Omega$  of the function
$\rho_t$ parametrizing $\Gamma_t$ as a normal graph over $\Gamma_0$.
 Section 3: Following the
argument in \cite{BaconneauLunardi}, we introduce a transformation
of the dependent variables:
$$u^I(x,t),\tilde{\rho}(x,t)\longrightarrow U^I(x,t),\quad x\in D_0^I,$$
 which leads to a
system on $\Omega\times [0,T]$ with three unknown functions and
three conjugation conditions on $\Gamma_0$. Section 4: We describe
the fixed-point scheme, the associated linear problem and the
required H\"{o}lder space estimates. This concludes the proof,
except for verifying that the linear system satisfies the
`complementarity conditions', which is done in Section
5.\vspace{.2cm}

\emph{Acknowledgments.} This research was started during a semester
visit to the Max-Planck Institute for Gravitational Physics
(Potsdam-Golm); I am grateful to the Geometric Analysis group for
the invitation and to the Max-Planck Society for support. An early
version of this result was announced at a seminar at the FU Berlin
(June 2007), and the current version at the conference Nonlinear
PDEs at IMPA (August 2008). Finally, it is a pleasure to thank
N.Alikakos (University of Athens) for originally suggesting the
problem and for his interest in the results.

\vspace{.5cm} \textbf{1. Compatibility conditions.}\vspace{.2cm}

We wish to find solutions that converge to the initial data in a
strong sense, as $t\rightarrow 0_+$. It is well-known from parabolic
theory that, for such solutions to exist, the initial data have to
satisfy compatibility conditions; in the present case, of `order 0'
and `order 1'. Before going any further, one should check that the
geometric problem is compatible with a graph formulation, in the
sense that one can find reasonably large families of initial
configurations satisfying these compatibility requirements. In this
section we state the conditions, analyze their geometry and examine
their internal consistency.\vspace{.2cm}

`Order 0' compatibility is just the requirement that the initial
data satisfy the conjugation and boundary conditions (i.e., that the
$w^I_0$ define a configuration).\vspace{.2cm}

On $\Gamma_0$ this corresponds to two zero-order conditions on the
$w_0^I$ (pointwise matching) and two first-order conditions (angle).
Assuming the ${w_0^I}_{|\Gamma_0}$ are given (and hence the curve
$\Lambda_0$), the conditions on the first-order data (the tangent
planes to $\Sigma^I(0)$ along $\Lambda(0)$) are clearly compatible.
(Construction: at each point, on the plane orthogonal to the tangent
line of $\Lambda_0$- which is determined by zero-order data- we can
arbitrarily pick one unit vector, say $N^1(0)$, and then the other
two are uniquely determined; so we have one real degree of freedom.)
\vspace{.3cm}

`Order one' compatibility conditions along $\Gamma_0$ are obtained
by differentiating the incidence conditions in $t$ and using the
equation, bearing in mind that $\Gamma$ is time-dependent.

Parametrizing $\Gamma_t$ by $a(s,t)\in {\mathbb R}^2$ ($s\in
{\mathbb R}/2\pi \mathbb{Z}$ is a fixed periodic variable, not arc
length), we compute (using the equation of motion) the velocity
vector at points of $\Lambda_t$, $V^I=\frac{\partial}{\partial
t}(G^I(a(s,t),t))$:
$$V^I=[\dot{a},Dw^I\cdot \dot{a}+\partial_tw^I]=[\dot{a},D
w^I\cdot \dot{a}+v^I H^I],$$ where $\dot{a}=\partial_ta$.

The compatibility condition is then the statement that, at $t=0$
(and, \emph{a fortiori}, for all $t$):
$$V^1=V^2=V^3,\quad y\in \Gamma_t.$$
In principle, this corresponds to up to six independent scalar
conditions. To analyze the situation, observe that at points of
$\Lambda_t$ we have three natural orthonormal frames
$\{E,N^I,T^I\}$, where $E$ is the unit tangent vector to $\Lambda_t$
(orienting $\Gamma_t$ counterclockwise) and we define $T^I=E\wedge
N^I$ (vector product), a tangent vector to $\Sigma^I$; in
particular, $T^1+T^2=T^3$. Noting that $\langle V^I,N^I\rangle=H^I$,
we may write:
$$V^I=\lambda^I E+\mu^IT^I+H^IN^I.$$
Six independent scalar conditions equivalent to equality of the
$V^I$ are then obtained as follows:

(1.) Add $\langle V^1,N^1\rangle=\langle V^3,N^1\rangle$ and
$\langle V^2,N^2\rangle=\langle V^3,N^2\rangle$.

(2.)  Add $\langle V^1,T^1\rangle=\langle V^3,T^1\rangle$ and
$\langle V^2,T^2\rangle=\langle V^3,T^2\rangle$.

(3.) Take the difference $\langle V^1,N^1\rangle=\langle
V^3,N^1\rangle$ minus $\langle V^2,N^2\rangle=\langle
V^3,N^2\rangle$.

(4.) Take the difference $\langle V^1,T^1\rangle=\langle
V^3,T^1\rangle$ minus $\langle V^2,T^2\rangle=\langle
V^3,T^2\rangle$.

(5./6.) $\langle V^1,E\rangle=\langle V^2,E\rangle=\langle
V^3,E\rangle$.

As we proceed to write these relations in a more explicit form, we
will preserve their labels (1)...(6). We remark that conditions
(1)...(4) already occur for curves, while (5) and (6) are new for
surfaces.

In terms of the coefficients $\mu^I$ and $\lambda^I$, we easily
obtain (using $N^1+N^2=N^3$ and $T^1+T^2=T^3$): \vspace{.2cm}

(1)$H^1+H^2=H^3$.

(2)$\mu^1+\mu^2=\mu^3$.

(3)$H^1-H^2=\sqrt{3} \mu^3$.

(Here we used $\langle T^3,N^1-N^2\rangle=\sqrt{3}$ and $\langle
N^3,N^1-N^2\rangle=0$.)

(4)$\mu^1-\mu^2=-\sqrt{3} H^3$.

(5/6) $\lambda^1=\lambda^3$ and $\lambda^2=\lambda^3$.

To make these relations explicit, one has to compute the $T^I$,
$\lambda^I$ and $\mu^I$, which is elementary. Denote by $\tau$ the
unit tangent vector to $\Gamma_t$, and set
$v_E=\sqrt{1+(d_{\tau}w^I)^2}$ (independent of $I$ on $\Gamma_t$)
and $(D w^I)^{\perp}=(-\partial_2w^I,\partial_1w^I)$. The results
are:
$$T^I=-\frac 1{v^Iv_E}[n-(d_{\tau}w^I)(Dw^I)^{\perp},d_n
w^I],$$
$$\mu^I=-\frac {1}{v_E}((\dot{a}\cdot
n)v^I+(d_n w^I)H^I),$$
$$\lambda^I=\frac 1{v_E}(\dot{a}\cdot
{\tau}+(d_{\tau}w^I)(\dot{a}\cdot D w^I+v^IH^I)).$$

With this we can rewrite (1)-(4) as:\vspace{.2cm}

(1)$H^1+H^2=H^3$.

(2)$(\dot{a}\cdot n)(v^1+v^2)+(d_nw^1)H^1+(d_nw^2)H^2=(\dot{a}\cdot
n)v^3+(d_n w^3)H^3$.

 (3) $H^1-H^2=(\sqrt{3}/v_E)((\dot{a}\cdot
 n)v^3+(d_nw^3)H^3)$.

 (4) $(\dot{a}\cdot n)(v^1-v^2)+(d_n
 w^1)H^1-(d_n w^2)H^2=-\sqrt{3}v_EH^3$.
 \vspace{.2cm}

 When rewriting conditions (5/6), we observe that terms involving
 the tangential component $\dot{a}\cdot \tau$ of the moving boundary
 velocity cancel, and there is a common factor $d_{\tau}w\equiv d_{\tau}w^I$ (independent of $I$). We have:
 \vspace{.2cm}

 (5/6) $(d_{\tau} w)((\dot{a}\cdot
 n)d_nw^1+v^1H^1)=(d_{\tau} w)((\dot{a}\cdot
 n)d_nw^2+v^2H^2)$
 $$=(d_{\tau} w)((\dot{a}\cdot
 n)d_nw^3+v^3H^3).$$

 \vspace{.2cm}

 Thus (5) and (6) occur only at points where $d_Tw\neq 0$
 (that is, where $\Lambda_t$ is not `horizontal'). Regarding the zero
 and first order data of the configuration as given, (1)-(6) give
 six homogeneous linear relations for the four `second-order
 quantities' $\dot{a}\cdot n$ and $H^I$, so it is conceivable that
 only the zero solution might exist. That this is not so is part of the
 conclusion of the following proposition.
\vspace{.5cm}

\textbf{Proposition 1.1.} \emph{The homogeneous system of linear
relations (1)-(6) has rank two, and is generated by (1) and (3).}
\vspace{.3cm}

This has the following interpretation: given zero and first-order
data, the only constraint on the second-order data of an initial
configuration of surfaces satisfying `order one compatibility at
$t=0$' is condition (1) on the mean curvatures at the junction
curve, while (3) specifies how the mean curvatures determine the
normal velocity of the interface $\Gamma_t$, in particular at $t=0$
(since the normal velocity $\dot{a}\cdot n$ is independent of how we
parametrize $\Gamma_t$, we denote it by $\dot{\Gamma}_n$):

$$\dot{\Gamma}_n=\frac{v_E}{v^3\sqrt 3}(H^1-H^2)-\frac{\partial_n w^3}{v^3}H^3.$$

Geometrically in $\mathbb{R}^3$, condition (3) says that the
component $V^{nor}$ (normal to the junction $\Lambda_t$) of the
velocity vector, which is independent of how one parametrizes
$\Lambda_t$, is determined by the $H^I$ via:
$$V^{nor}=H^3N^3+\mu^3T^3=H^3N^3+\frac 1{\sqrt{3}}(H^1-H^2)T^3.$$

\textbf{Proof.} This is computational and elementary, so we just
describe the main steps. One just has to observe that (as for
curves), the constraints on first derivatives arising from writing
$N^1+N^2=N^3$ componentwise can be `solved', in the sense of the
following elementary algebraic lemma: \vspace{.2cm}

\textbf{Lemma 1.2.} Given real numbers $\alpha>0$ and $a^1,a^2,a^3$,
let $v^I=\sqrt{\alpha^2 +(a^I)^2}$; suppose these numbers satisfy
the following relations:
$$\frac 1{v^1}+\frac 1{v^2}=\frac
1{v^3},\quad \frac{a^1}{v^1}+\frac{a^2}{v^2}=\frac{a^3}{v^3}.$$ Then
(assuming $a^2>a^1$ to normalize the labeling):
$$\frac 1{v^1}=\frac
1{v^3}(\frac 12+\frac{\sqrt 3}2\frac{a^3}{\alpha}),\quad \frac
1{v^2}=\frac 1{v^3}(\frac 12-\frac{\sqrt 3}2\frac{a^3}{\alpha}) ,$$
$$\frac{a^1}{v^1}=\frac 12\frac {a^3}{v^3}-\frac{\sqrt 3}2\frac
{\alpha}{v^3},\quad \frac{a^2}{v^2}=\frac 12\frac
{a^3}{v^3}+\frac{\sqrt 3}2\frac {\alpha}{v^3}.$$

\emph{Remark.} This lemma has the following equivalent geometric
formulation: let $\omega^1,\omega^2,\omega^3$ be vectors in
$\mathbb{R}^2$ of the same length, satisfying
$\omega^1+\omega^2=\omega^3$. Then, denoting by $R_{\theta}$ the
counterclockwise rotation operator in $\mathbb{R}^2$ (by $\theta$
radians), we have:
$$\omega^1=R_{\pi/3}[\omega^3],\quad \omega^2=R_{5\pi/3}[\omega^3].$$
\vspace{.2cm}

Since $n$ is the inner unit normal of $D^1_t=D^2_t$, the condition
$a^2>a^1$ corresponds to labeling the surfaces so that $\Sigma^2$ is
below $\Sigma^1$ near $\Lambda_t$.\vspace{.2cm}

To prove proposition 1.1, we substitute the values of $v^I$ and
$d_nw^I/v^I$ given by lemma 1.2 (for $I=1,2$, in terms of the same
quantities for $I=3$) into relations (2)-(6). Using relation (1), we
find that each of (2) and (4) is equivalent to (3). Where
$d_{\tau}w\neq 0$, adding and subtracting the relations obtained
from (5) and (6) with this substitution, we recover respectively (1)
and (3), completing the proof.

\vspace{.2cm}

There is no `order 1' compatibility condition to consider at the
fixed boundary $\partial \Omega$, since the Neumann condition
$d_{\omega}w^3=0$ is first-order. \vspace{.5cm}

\textbf{2. An equivalent system on a fixed domain.} \vspace{.2cm}

In this section we describe new dependent variables
$u^I(x,t),\tilde{\rho}(x,t)$ defined on fixed cylindrical domains
$D_0^I\times [0,T]$, $\Omega \times [0,T]$ (resp.), and a parabolic
system with conjugation conditions, equivalent to the original
system for $w^I(y,t)$ on $D^I_t, t\in [0,T]$.\vspace{.2cm}

\textbf{2.1. An extension operator.} $\Gamma_0\subset \Omega\subset
\mathbb{R}^2$ is a simple closed curve of class $C^{3+\alpha}$,
bounding domains $D_0^1=D_0^2$ on the inside and
$D_0^3=\Omega-\bar{D}_0^1$; denote by $\nu$ the inner unit normal
vector field of $\Gamma_0$ (of class $C^{2+\alpha}$) and by $r:{\cal
N}_0\rightarrow \mathbb{R}$ the oriented distance to $\Gamma_0$,
where ${\cal N}_0$ is a fixed tubular neighborhood of $\Gamma_0$. We
have $r\in C^{3+\alpha}({\cal N}_0)$, and the natural extension
$\nu=Dr$ of $\nu$ to ${\cal N}_0$. With $\bar{r}:=\sup_{{\cal
N}_0}r$, consider the smaller tubular neighborhood:
$${\cal N}:=\{x\in {\cal N}_0;|r(x)|<r_0:=\min\{\bar{r}/2,1/2\}\}.$$
 Denoting by $\hat{\zeta}\in C^{\infty}(\mathbb{R};[0,1])$ a
fixed function supported in $[-r_0,r_0]$, equal to $1$ in
$[-r_0/2,r_0/2]$ with $|\hat{\zeta}'|\leq 3/r_0$ everywhere, we see
that $\zeta(x)$, defined as $\hat{\zeta}\circ r$ in $\cal N$ and
identically zero outside of $\cal N$, is in
$C^{3+\alpha}(\mathbb{R}^2)$ and coincides with $r$ in ${\cal
N}_1:=r^{-1}(-r_0/2,r_0/2)\subset {\cal N}$.\vspace{.2cm}

Using $\zeta$ we extend $\nu=Dr$ from ${\cal N}_1$ to
$\mathbb{R}^2$, by setting $\bar{\nu}:=\zeta\nu\in
C^{2+\alpha}(\mathbb{R}^2; \mathbb{R}^2)$. $\bar{\nu}$ is supported
in $\cal N$ and satisfies $|\bar{\nu}|\leq 1$ in $\mathbb{R}^2$.
Note that $D_{\nu}\nu\equiv 0$ in ${\cal N}_1$, while:
$$D_{\bar{\nu}}\bar{\nu}=\zeta D_{\nu}(\zeta \nu)=\zeta
(d_{\nu}\zeta) \nu=(d_{\nu}\zeta) \bar{\nu},$$ which is supported in
${\cal N}-{\cal N}_1$.\vspace{.2cm}

We will need a fixed extension operator
$C^{k+\alpha}(\Gamma_0)\rightarrow C^{k+\alpha}(\mathbb{R}^2)$,
defined in a standard way using $\zeta$ (for $k\geq 2$). Given
$\rho\in C^{k+\alpha}(\Gamma_0)$, extend $\rho$ to ${\cal N}_0$ by
projection along normal line segments (so $d_{\nu}\rho\equiv 0$ in
${\cal N}_0$), then set:
$$\tilde{\rho}:=\left \{ \begin{array}{l}\zeta
\rho,\quad x\in {\cal N}_0\\0,\quad x\in {\mathbb R}^2-{\cal
N}\end{array}\right .$$ It is clear that $\tilde{\rho}$ is supported
in $\cal N$, and that
$||\tilde{\rho}||_{C^{k+\alpha}(\mathbb{R}^2)}\leq c_{\zeta}
||\rho||_{C^{k+\alpha}(\Gamma_0)}$, where $c_{\zeta}$ depends only
on $\zeta$ and on $||r||_{C^{3+\alpha}(\bar{\cal N}_0)}$. Since the
extension operator is $t$-independent, the same estimate holds for
parabolic H\"{o}lder norms when $\rho$ depends on $t$. Note also:
$$D_{\bar{\nu}}(\tilde{\rho}\bar{\nu})=\zeta
\nu(\tilde{\rho})\bar{\nu}+\tilde{\rho}D_{\bar{\nu}}\bar{\nu}$$
$$=\zeta(\nu(\zeta)\rho+\zeta\nu(\rho))\bar{\nu}+\tilde{\rho}(d_{\nu}\zeta)\bar{\nu}$$
$$=\zeta((d_{\nu}\zeta)\rho+0)\bar{\nu}+\tilde{\rho}(d_{\nu}\zeta)\bar{\nu}$$
$$=2(d_{\nu}\zeta)\tilde{\rho}\bar{\nu},$$
supported in ${\cal N}-{\cal N}_1$.

If $U$ is already defined on a set containing $\Gamma_0$, we denote
by ${\cal E}[U]$ the restriction-extension operator:
$${\cal E}[U]:=(U_{|\Gamma_0})^{\tilde{}}.$$

\vspace{.3cm}

\textbf{2.2. Diffeomorphism.} For $t>0$ small enough, the interface
$\Gamma_t$ will be contained in ${\cal N}_1$, and we may parametrize
it as a normal graph over $\Gamma_0$:
$$\Gamma_t=\{x+\rho(x,t)\nu(x); x\in \Gamma_0\},\mbox{ for some
}\rho:\Gamma_0\times [0,T]\rightarrow \mathbb{R}.$$ Recall the
subscript $t$ always denotes `function (or map, or set) at time
$t$', so $\rho_t(x):=\rho(x,t)$. We use the extension operator $\cal
E$ to define a diffeomorphism $\varphi_t$ of $\mathbb{R}^2$ (or of
$\Omega$), for $t\in [0,T]$:
$$\varphi_t(x)=x+\tilde{\rho}_t(x)\bar{\nu}(x),\quad \varphi_0=id.$$
For $t\in [0,T]$ small enough, $\varphi_t$ is a diffeomorphism (of
class $C^{2+\alpha}$), equal to the identity in ${\mathbb R}^2-{\cal
N}$, and mapping:
$$\Omega\rightarrow \Omega, \quad\Gamma_0\rightarrow \Gamma_t,\quad D_0^I\rightarrow D^I_t,\quad
I=1,2,3,$$ diffeomorphically in each case. If $\rho\in
C^{2+\alpha,1+\alpha/2}(\Gamma_0\times [0,T])$, $\varphi$ is in the
same H\"{o}lder class as a function of $(x,t)$. The differential of
$\varphi_t$ is given by:
$$D\varphi_t=\mathbb{I}+(D\tilde{\rho}_t)\bar{\nu}+\tilde{\rho}_tD\bar{\nu};$$
in particular, in light of the calculation in 2.1, $D\varphi_t$ maps
$\bar{\nu}$ via:
$$D\varphi_t[\bar{\nu}]=\bar{\nu}+D_{\bar
\nu}(\tilde{\rho}_t\bar{\nu})=[1+2(d_{\nu}\zeta)\tilde{\rho}_t]\bar{\nu},$$
which implies:
$$({D\varphi_t})^{-1}[\bar{\nu}]=\frac
1{1+2(d_{\nu}\zeta)\tilde{\rho}_t}\bar{\nu}:=\bar{\zeta}(\tilde{\rho})\bar{\nu}.$$
Here $\bar{\zeta}(\tilde{\rho})\in C^{3+\alpha}(\mathbb{R}^2)$, and
is identically $1$ in $\bar{\cal N}_1\cap ({\mathbb R}^2-{\cal N})$,
with modulus bounded above by $1$ everywhere (provided only we have
$|\tilde{\rho}_t|\leq r_0/12$ for $t\in [0,T]$). In particular, on
$\Gamma_0$:
$${D\varphi_t}_{|\Gamma_0}:\nu(x)\mapsto \nu(x+\rho_t(x)\nu(x))$$
(recall $\nu:=Dr$ on ${\cal N}_0$). Of course,
$\nu(x+\rho_t(x)\nu(x))$ is in general not normal to
$\Gamma_t$.\vspace{.3cm}

\textbf{2.3 Equations over $D_0^I$.} Assume $w_t^I:D^I_t\rightarrow
\mathbb{R}$ are solutions of $(MCM)_w$ for $t\in [0,T]$. Define, for
$I=1,2,3$:
$$u^I(x,t)=w^I(\varphi_t(x),t),x\in D_0^I;\quad
w^I(y,t)=u^I(\psi_t(y),t),y\in D_t^I$$ (where
$\psi_t:=(\varphi_t)^{-1}$). A standard calculation yields (omitting
the superscript $I$ sometimes, and using the Einstein summation
convention throughout):
$$0=\partial_tw-g^{ij}(Dw)D^2_{i,j}w=\partial_tu-h^{ab}D^2_{a,b}u+(\partial_cu)
(\partial_t\psi^c-g^{ij}D^2_{i,j}\psi^c),$$ or equivalently:
$$\partial_tu-h^{ab}D^2_{a,b}u-(\partial_cu)(\partial_i\psi^c)
(\partial_t\varphi^i-h^{ab}\partial^2_{a,b}\varphi^i)=0,$$ where
$h_{ab}=\partial_a\varphi^i\partial_b\varphi^jg_{ij}$ is the
pullback metric (of $g$ under $\varphi$). Since the parametrization
over $D_0^I$ of $\Sigma^I$ is:
$$F_t^I(x)=G_t^I(\varphi_t(x))=[\varphi_t(x),w_t^I(\varphi_t(x))]=[\varphi_t(x),u_t^I(x)]
\in {\mathbb R}^2\times \mathbb{R},$$ we have for $x\in D_0^I$:
$$h_{ab}^I=\langle \partial_a
F_t^I,\partial_bF_t^I\rangle=\partial_a\varphi_t\cdot
\partial_b\varphi_t+(\partial_au_t^I)(\partial_bu_t^I)$$
$$=\delta_{ab}+\partial_a(\tilde{\rho}_t\bar{\nu})\cdot
\partial_b(\tilde{\rho}_t\bar{\nu})+(\partial_au_t^I)(\partial_bu_t^I)$$
$$:=\delta_{ab}+h'_{ab}(\tilde{\rho}_t,D\tilde{\rho}_t,Du_t^I).$$
Define the operator on functions (or maps) $f(x,t)$ on $D_0^I\times
[0,T]$:
$$L_{h^I}[f]:=\partial_tf-tr_{h^I}D^2f.$$
The equation on $D_0^I$ is:
$$L_{h^I}[u^I]-Du^I\cdot (D\varphi_t)^{-1}[L_{h^I}\varphi]=0.$$
We may express $L_{h}[\varphi]$ in terms of $\tilde{\rho}$ and
$\bar{\nu}$:
$$L_h[\varphi]=L_h[\tilde{\rho}]\bar{\nu}-\tilde{\rho}tr_hD^2\bar{\nu}-2\langle
D\tilde{\rho},D\bar{\nu}\rangle_h.$$ Using $(D\varphi)^{-1}[\bar
\nu]=\bar{\zeta}(\tilde{\rho})\bar{\nu}$, we find:
$$(D\varphi)^{-1}[L_h\varphi]=\bar{\zeta}(\tilde{\rho})L_h[\tilde{\rho}]\bar{\nu}
-(D\varphi)^{-1}[\tilde{\rho}tr_hD^2\bar{\nu}+2\langle
D\tilde{\rho},D\bar{\nu}\rangle_h].$$ Hence the system in the
variables $(u^I,\tilde{\rho})$ is:
$$L_{h^I}[u^I]-\bar{\zeta}(\tilde{\rho})(d_{\bar{\nu}}u^I)L_{h^I}[\tilde{\rho}]
+Du^I\cdot (D\varphi_t)^{-1}[\tilde{\rho}tr_{h^I}D^2\bar
{\nu}+2\langle D\tilde{\rho},D\bar{\nu}\rangle_{h^I}]=0,\quad
(PDE)_{(u^I,\rho)}$$ where:
$$h_{ab}^I=\delta_{ab}+h'_{ab}(\tilde{\rho},D\tilde{\rho},Du^I),\quad
D\varphi_t=\mathbb{I}+D(\tilde{\rho}\bar{\nu}),$$
$$\langle
D\tilde{\rho},D\bar{\nu}\rangle_{h^I}=(h^I)^{ab}(\partial_a\tilde{\rho})(\partial_b\bar{\nu})\in
\mathbb{R}^2.$$\vspace{.3cm}

\textbf{2.4. Boundary conditions at $\Gamma_0$.} (More precisely,
these are `conjugation conditions', since $\Gamma_0$ is an `internal
boundary' in $\Omega$). We begin by computing how
${D\varphi_t}_{|\Gamma_0}$ acts on the orthonormal frame
$\{\nu,\tau_0\}$, $\tau_0:=-\nu^{\perp}$ (where $\perp$ denotes
counterclockwise $\pi/2$ rotation in $\mathbb{R}^2$). Parametrizing
$\Gamma_0$ by arc length $s$, we have:
$$\Gamma_0'(s)=\tau_0,
\quad \nu'=-k_0\tau_0$$ (with $k_0$ the curvature of $\Gamma_0$).
Composing with $\varphi_t$, we obtain the parametrization of
$\Gamma_t$: $\Gamma_t(s)=\Gamma_0(s)+\rho_t(\Gamma_0(s))\nu(s),$
with tangent and inner unit normal vectors:
$$\Gamma_t'(s)=(1-k_0\rho)\tau_0+(d_{\tau_0}\rho)\nu,\quad
n=\frac{(1-k_0\rho)\nu-(d_{\tau_0}\rho)\tau_0}{[(1-k_0\rho)^2+d_{\tau_0}\rho)^2]^{1/2}}.$$
Using $D_{\nu}\nu=0,d_{\nu}\rho=0,D_{\tau_0}\nu=\nu'$ on $\Gamma_0$,
we find for an arbitrary $v\in \mathbb{R}^2$ and $x\in \Gamma_0$:
$$D\varphi_t(x)[v]=v+(v\cdot
\tau_0)(d_{\tau_0}\rho)\nu-(v\cdot \tau_0)(k_0\rho)\tau_0,$$ in
particular verifying again:
$$D\varphi_t[\nu]=\nu, \quad
D\varphi_t[\tau_0]=(1-k_0\rho)\tau_0+(d_{\tau_0}\rho)\nu=\Gamma_t'(s).$$
We need vectors mapping to $n$ and to $\tau:=-n^{\perp}$ under
$D\varphi_t$. It is easy to see that, defining:
$$\hat{\tau_0}:=[(1-k_0\rho)^2+(d_{\tau_0}\rho)^2]^{-1/2}\tau_0:=a_{11}\tau_0,$$
$$\mu=(1-k_0\rho)^{-1}\{[(1-k_0\rho)^2+(d_{\tau_0}\rho)^2]^{1/2}\nu-(d_{\tau_0}\rho)\hat{\tau_0}\}$$
$$:=a_{21}\tau_0+a_{22}\nu,$$ we have for these vector fields
$\mu=\mu(\rho,d_{\tau_0}\rho),\hat{\tau_0}=\hat{\tau_0}(\rho,d_{\tau_0}\rho)$:
$$D\varphi_t(x)[\mu]=n,\quad D\varphi_t(x)[\hat{\tau_0}]=\tau,\quad
x\in \Gamma_0.$$ (The expressions with the $a_{ij}$ are used in
section 4.) Hence, for $x\in \Gamma_0$:
$$d_nw^I(\varphi_t(x))=d_{\mu}u^I(x);\quad
d_{\tau}w^I(\varphi_t(x))=d_{\hat{\tau_0}}u^I(x),$$
$$v^I=[1+(d_{\tau}w^I)^2+(d_nw^I)^2]^{1/2}=[1+(d_{\hat{\tau_0}}u^I)^2+(d_{\mu}u^I)^2]^{1/2}$$
$$:={\cal G}(\rho,d_{\tau_0}\rho,d_{\tau_0}u^I,d_{\nu}u^I).$$
(\emph{Remark:} It is easy to see that
$\mu-\nu=-(d_{\tau_0}\rho)\tau_0+O(\rho^2+|d_{\tau_0}\rho|^2)$; at
$t=0$, $\mu\equiv \nu$ on $\gamma_0$, since $\rho\equiv
0$.)\vspace{.2cm}

We now state the conjugation conditions in terms of $u^I,\rho$:
$$\left \{\begin{array}{lr}
u^1(x,t)=u^2(x,t)=u^3(x,t) \quad &(BC^0)_{(u^I,\rho)}\\
\frac 1{v^1}+\frac 1{v^2}=\frac 1{v^3}\quad &(BC^1)_{(u^I,\rho)}\\
\frac{d_{\nu}u^1}{v^1}+\frac{d_{\nu}u^2}{v^2}=\frac{d_{\nu}u^3}{v^3}\quad
&(BC^2)_{(u^I,\rho)} \end{array}\right .$$
 In the computation of the
second angle condition, we use the fact that $d_{\hat{\tau_0}}u^I$
is independent of $I$ on $\Gamma_0$ (which follows from the matching
condition). Here $v^I={\cal
G}(\rho,d_{\tau_0}\rho,d_{\tau_0}u^I,d_{\nu}u^I)$ has the expression
given above. \vspace{.3cm}

Since the diffeomorphism $\varphi_t$ is the identity on $\partial
\Omega$, we have the additional (Neumann) boundary condition for
$u^3$:
$$\partial_{\omega}u^3=0\mbox{ on }\partial\Omega.$$
(Recall $\partial D^3_0=\Gamma_0\sqcup \partial \Omega$.)

\vspace{.5cm} \textbf{3. The transformed system.}

The system described in the preceding section includes 4 scalar
conjugation conditions, but only three evolution equations for the
four unknowns $u^I,\rho$; there is no explicit evolution equation
for $\rho$ (or its extension $\tilde{\rho}$). Yet (assuming the
evolution of triple junctions is well-posed as a PDE) it is clear
geometrically that the evolution of the $u^I$ determines that of
$\rho$. In this section, following the technique introduced in the
paper \cite{BaconneauLunardi} we show it is possible to define new
dependent variables $U^I$ on $D_0^I\times [0,T]$ in such a way that
$\tilde{\rho}$ is recovered from the $U^I$ via the extension
operator $\cal E$. The drawback (as in \cite{BaconneauLunardi}) is
having to introduce `non-local highest-order terms' in the resulting
equation for the $U^I$. \vspace{.2cm}

Following \cite{BaconneauLunardi}, define $U^I:D^I_0\times
[0,T]\rightarrow \mathbb{R}$ via:
$$u^I(x,t)=u_0^I(x)+(d_{\bar
{\nu}}u_0^I)(x)\tilde{\rho}(x,t)+U^I(x,t).$$ (Note $U^I_{|t=0}\equiv
0$.) The matching conditions $(BC^0)_{(u^I,\rho)}$ on $\Gamma_0$
imply:
$$\rho d_{\nu}u_0^1+U^1=\rho d_{\nu}u_0^2+U^2=\rho
d_{\nu}u_0^3+U^3\mbox{ on }\Gamma_0.$$ Equivalently, we have the
following equalities on $\Gamma_0$:
$$\quad \rho_{|\Gamma_0}=\frac{U^2-U^1}{d_{\nu}u_0^1-d_{\nu}u_0^2}=
\frac{U^3-U^2}{d_{\nu}u_0^2-d_{\nu}u_0^3}=\frac{U^3-U^1}{d_{\nu}u_0^1-d_{\nu}u_0^3}.$$
The last equality may be written in the form:
$$\frac{U^3-U^2}{U^3-U^1}=\frac{d_{\nu}u_0^2-d_{\nu}u_0^3}{d_{\nu}u_0^1-d_{\nu}u_0^3}
=-\frac{v_0^2}{v_0^1},$$ which we regard as a new form of the
matching condition. We used the fact that the angle conditions at
$t=0$ imply:
$$\frac 1{v_0^1}(d_{\nu}u_0^1-d_{\nu}u_0^3)=\frac
1{v_0^2}(d_{\nu}u_0^3-d_{\nu}u_0^2).$$ This new matching condition
may also be written in the form:
$$\frac{U^3-U^2}{v_0^2}+\frac{U^3-U^1}{v_0^1}=0,$$
or equivalently (using the first angle condition at $t=0$):
$$\frac{U^1}{v_0^1}+\frac{U^2}{v_0^2}=\frac{U^3}{v_0^3}\mbox{ on
}\Gamma_0\quad (BC^0)_{U^I},$$ which we adopt as the matching
condition for the $U^I$ system.\vspace{.2cm}

The main point of the method is that the relation:
$$\rho=\delta_0(U^2-U^1)\mbox{ on }\Gamma_0,\quad
\delta_0:=(d_{\nu}u_0^1-d_{\nu}u_0^2)^{-1}$$ allows us to express
$\tilde{\rho}$ in terms of ${\cal E}[U^2-U^1]$.\vspace{.2cm}

To explain how this is done, we take a moment to examine the
extension operator applied to the product of two functions, $f,g$
defined on $\Gamma_0$. Denote by the subscript $rad$ the extension
to ${\cal N}_0$ constant along normals of a function defined on
$\Gamma_0$. By definition, on ${\cal N}_0$ we have:
$$(fg)^{\tilde{}}=\zeta(fg)_{rad}=(f)_{rad}\tilde{g}$$ (this is
true on $\cal N$, and both sides vanish on ${\cal N}_0-{\cal N}$).
Let $\chi$ be a smooth function in $\mathbb{R}^2$, equal to 1 on
$\cal N$ and vanishing outside ${\cal N}_0$.  Given $f$ on
$\Gamma_0$, define $\hat{f}$ on $\mathbb{R}^2$ by:
$$\hat{f}=\chi f_{rad}\mbox{ on }{\cal N}_0,\quad \hat{f}=0\mbox{ on
}\mathbb{R}^2-{\cal N}_0.$$ Then it is easy to check that:
$$(fg)^{\tilde{}}=\hat{f}\tilde{g}\mbox{ on }\mathbb{R}^2.$$
\vspace{.2cm}

Thus we may recover $\tilde{\rho}$ from $U^2-U^1$ as follows:
$$\tilde{\rho}=\hat{\delta_0}{\cal E}[U^2-U^1].$$

We now express the  angle boundary conditions in terms of the $U^I$.
Since:
$$u^I=u_0^I+\delta_0(U^2-U^1)d_{\nu}u_0^I+U^I\mbox{ on }\Gamma_0,$$
we have on $\Gamma_0$:
$$d_{\tau_0}u^I=d_{\tau_0}u_0^I+d_{\tau_0}[\delta_0(U^2-U^1)d_{\nu}u_0^I]+d_{\tau_0}U^I$$
$$:=A_{tan}(U^2-U^1,d_{\tau_0}(U^2-U^1))+d_{\tau_0}U^I$$
($A_{tan}$ for `affine', as a function of the indicated arguments.)
For the normal derivatives, we use the expression defining $U^I$ on
$D_0^I$, and recall that $d_{\nu}\tilde{\rho}=0$ on $\Gamma_0$;
hence (using also $D_{\nu}\nu=0$):
$$d_{\nu}u^I_{|\Gamma_0}=d_{\nu}u_0^I+\delta_0(U^2-U^1)D^2u_0^I(\nu,\nu)+d_{\nu}U^I$$
$$:=A_{nor}(U^2-U^1)+d_{\nu}U^I.$$
In section 2 we found functional expressions $v^I_{|\Gamma_0}={\cal
G}({\rho},d_{\tau_0}\rho, d_{\tau_0}u^I,d_{\nu}u^I)$. If in the
expressions defining the vector fields $\hat{\tau}_0$ and $\mu$ on
$\Gamma_0$ we replace $\rho$ and $d_{\tau_0}\rho$ by
$\delta_0(U^2-U^1)$ and $d_{\tau_0}[\delta_0(U^2-U^1)]$ (resp.) and
use the expressions just found for $d_{\tau_0}u^I$ and $d_{\nu}u^I$,
we obtain:
$$v^I_{|\Gamma_0}={\cal G}(\delta_0(U^2-U^1), d_{\tau_0}[\delta_0(U^2-U^1)],
 A_{tan}(U^2-U^1,d_{\tau_0}(U^2-U^1))+d_{\tau_0}U^I,
A_{nor}(U^2-U^1)+d_{\nu}U^I).$$ \vspace{.2cm}

From now on we adopt the notational convention for the summation of
quantities depending on a superscript $I=1,2,3$:
$${\sum_I}'a^I:=a^1+a^2-a^3.$$\vspace{.2cm}

It follows from the above (using also the boundary conditions for
the $u_0^I$) that the boundary conditions take the form:
$$\left \{ \begin{array}{lr}
\sum_I'\frac{U^I}{v_0^I}=0\quad &(BC^0)_{U^I}\\
\sum_I'\frac 1{v^I}=0\quad &(BC^1)_{U^I}\\
\sum_I'
\frac{d_{\nu}U^I}{v^I}+\sum_I'\frac{d_{\nu}u_0^I}{v^I}+\delta_0(U^2-U^1)\sum_I'\frac{D^2u_0^I(\nu,\nu)}{v^I}=0\quad
&(BC^2)_{U^I}
\end{array}\right .$$

Computing the equation satisfied by the $U^I$ is straightforward.
First we write down the expression for $L_{h^I}[u^I]$, isolating the
lower-order terms:
$$L_{h^I}[u^I]=L_{h^I}[U^I]+(d_{\bar
{\nu}}u_0^I)L_{h^I}[\tilde{\rho}]+{\cal
L}^{(1)I}_{h^I}(\tilde{\rho},D\tilde{\rho}),$$ where:
$${\cal L}^{(1)I}_{h^I}(\tilde{\rho},D\tilde{\rho}):=-2\langle
D\tilde{\rho},D(d_{\bar{\nu}}u_0^I)\rangle_{h^I}-tr_{h^I}D^2u_0^I-\tilde{\rho}tr_{h^I}D^2(d_{\bar{\nu}}u_0^I).$$
Now using the equation $(PDE)_{(u^I,\rho)}$ from section 2 and the
definitions:
$${\cal
C}_{h^I}(\tilde{\rho},D\tilde{\rho}):=(D\varphi)^{-1}[\tilde{\rho}tr_{h^I}D^2\bar{\nu}+2\langle
D\tilde{\rho},D\bar{\nu}\rangle_{h^I}]\mbox{ with
}D\varphi=\mathbb{I}+D(\tilde{\rho}\bar{\nu}),$$
$${\cal
L}^{(2)I}_{h^I}(\tilde{\rho},D\tilde{\rho}):=D(u_0^I+\tilde{\rho}d_{\bar{\nu}}u_0^I)\cdot
{\cal C}_{h^I},$$ we find:
$$L_{h^I}[U^I]+[(d_{\bar{\nu}}u_0^I-\bar{\zeta}(\tilde{\rho})d_{\bar{\nu}}(u_0^I+\tilde{\rho}d_{\bar{\nu}}u_0^I+U^I)]
L_{h^I}[\tilde{\rho}]$$
$$+DU^I\cdot {\cal C}_{h^I}(\tilde{\rho},D\tilde{\rho})+{\cal L}^{(1)I}_{h^I}(\tilde{\rho},D\tilde{\rho})
+{\cal L}^{(2)I}_{h^I}(\tilde{\rho},D\tilde{\rho})=0.$$ (Note that
$(D\varphi)^{-1}$ depends on $D\tilde{\rho}$, and that the
first-order terms involving the metric $h^I$ depend on $DU^I$, in
addition to $\tilde{\rho}$ and $D\tilde{\rho}$.)\vspace{.2cm}

This may be regarded as a quasilinear system in $U^I$, if we agree
to replace $\tilde{\rho}$ by $\hat{\delta_0}{\cal E}[U^2-U^1]$ at
every occurrence of $\tilde{\rho}$ or $D\tilde{\rho}$. Due to the
presence of the restriction-extension operator ${\cal E}$, this
introduces `non-local terms', even to highest order.\vspace{.2cm}

To make the dependence on $U^I$ a bit more explicit, compute:
$$L_{h^I}[\tilde{\rho}]=L_{h^I}[\hat{\delta_0}{\cal E}[U^2-U^1]]$$
$$=\hat{\delta_0}L_{h^I}[{\cal E}[U^2-U^1]]-2\langle
D\hat{\delta_0},D({\cal E}[U^2-U^1])\rangle_{h^I}-{\cal
E}[U^2-U^1]tr_{h^I}D^2\hat{\delta_0}.$$

Thus we have the equation for $U^I$:
$$L_{h^I}[U^I]+{\cal
A}(\tilde{\rho},D\tilde{\rho},DU^I)_{|\tilde{\rho}=\hat{\delta_0}{\cal
E}[U^2-U^1]}L_{h^I}[{\cal E}[U^2-U^1]]$$
$$={\cal
F}_{U}^I(\tilde{\rho},D\tilde{\rho},{\cal E}[U^2-U^1],D({\cal
E}[U^2-U^1]))_{|\tilde{\rho}=\hat{\delta_0}{\cal E}[U^2-U^1]},$$
where:
$${\cal A}(\tilde{\rho},D\tilde{\rho},DU^I)=\hat{\delta_0}[(d_{\bar{\nu}}u_0^I
-\bar{\zeta}(\tilde{\rho})d_{\bar{\nu}}(u_0^I+\tilde{\rho}d_{\bar{\nu}}u_0^I+U^I)];$$
$${\cal F}^I_{U}=2\langle
D\hat{\delta_0},D({\cal E}[U^2-U^1])\rangle_{h^I}+{\cal
E}[U^2-U^1]tr_{h^I}D^2\hat{\delta_0}$$
$$-D(u_0^I+\tilde{\rho}d_{\bar{\nu}}u_0^I+U^I)\cdot
{\cal C}_{h^I}$$
$$+2\langle
D\tilde{\rho},D(d_{\bar{\nu}}u_0^I)\rangle_{h^I}+tr_{h^I}D^2u_0^I+\tilde{\rho}tr_{h^I}D^2(d_{\bar{\nu}}u_0^I),$$
and ${\cal C}_{h^I}$ was given above.\vspace{.2cm}

We record this in abbreviated form (changing the notation slightly
to exhibit the dependence of $h^I$ and $\cal A$ on $U$):
$$L_{h^I_U}[U^I]+{\cal A}_U^IL_{h^I_U}[{\cal E}[U^2-U^1]]
={\cal F}^I_{U}\quad (PDE)_{U^I}$$

The critical term is ${\cal A}^I_UL_{h^I_U}[{\cal E}[U^2-U^1]]$, due
to the `non-local dependence' on $U^I$ (i.e., it is not given by the
action of a differential operator on the $U^I$, so we may not quote
standard linear parabolic theory at this point). This term has the
important property:
$${\cal A}^I_{|t=0}\equiv 0\mbox{ on }D_0^I.$$
This is easy to see from the definition, since ${\cal A}^I$ is
$\hat{\delta_0}$ multiplied by the expression:
$$d_{\bar{\nu}}u_0^I
-\bar{\zeta}(\tilde{\rho})d_{\bar{\nu}}(u_0^I+\tilde{\rho}d_{\bar{\nu}}u_0^I+U^I).$$
At $t=0$ we have $U^I\equiv 0$, $\tilde{\rho}\equiv 0$ (since
$\rho\equiv 0$) and therefore $\bar{\zeta}(\tilde{\rho})\equiv 1$.
This makes it possible to treat this linear system as a small
perturbation of a standard linear parabolic system, as long as $T$
is small. This was also the general strategy adopted in
\cite{BaconneauLunardi}, but we actually deal with this term in a
different way (in a sense, less sophisticated; see Section 4).
$\tilde{\rho}_{|t=0}=0$ also implies ${\cal C}_{h^I}=0$ at
$t=0$.\vspace{.2cm}

Finally, we record the dependence of the metric $h^I$ on $DU^I$:
$$h_{ab}^I=\delta_{ab}+\partial_a(\tilde{\rho}\bar{\nu})\cdot
\partial_b(\tilde{\rho}\bar{\nu})$$
$$+(\partial_au_0^I+\partial_a(\tilde{\rho}d_{\bar{\nu}}u_0^I)+\partial_aU^I)
(\partial_bu_0^I+\partial_b(\tilde{\rho}d_{\bar{\nu}}u_0^I)+\partial_bU^I),$$
where we set $\tilde{\rho}=\hat{\delta_0}{\cal
E}[U^2-U^1]$.\vspace{.2cm}

Note also the two properties of the term ${\cal F}^I_{U}$: it
involves derivatives of $U^I$ up to first order, and equals
$tr_{h^I_0}D^2u_0^I=v_0^IH_0^I$ at $t=0$.

\vspace{.5cm} \textbf{4. The associated linear system.} In this
section we describe the (standard) fixed-point argument used to
prove local existence, based on the solution of an associated linear
system with `non-local terms'. \vspace{.2cm}

\textbf{4.1 Linearization of the conjugation conditions.}

The angle conjugation conditions $(BC)_{U^I}$ involve sums of
$1/v^I$, so we begin by considering the linearization ${\cal L}$ of
$1/v$ at points of $\Gamma_0$, where:
$$v_{|\Gamma_0}=\sqrt{1+(d_{\mu}u)^2+(d_{\hat{\tau_0}}u)^2}=v(\rho,d_{\tau_0}\rho,DU)$$ and we linearize at
$\rho\equiv 0, U\equiv 0$. We find:
$${\cal L}(\frac 1v)=-\frac 1{v_0^3}[(d_{\nu}u_0){\cal
L}(d_{\mu}u)+(d_{\tau_0}u_0){\cal L}(d_{\hat{\tau_0}}u)]$$ From the
expressions for $\hat{\tau_0}$ and $\mu$ given in section 2:
$$\hat{\tau}_0=a_{11}\tau_0,\quad \mu=a_{21}\tau_0+a_{22}\nu,\quad {a_{ij}}_{|\rho\equiv 0}=\delta_{ij},$$
we have:
$${\cal
L}(d_{\hat{\tau_0}}u)=d_{\tau_0}U+(d_{\nu}u_0)d_{\tau_0}\rho+{\cal
L}(a_{11}) d_{\tau_0}u_0+d^2u_0(\tau_0,\nu)\rho,$$
$${\cal L}(d_{\mu}u)=d_{\nu}U+{\cal L}(a_{21})d_{\tau_0}u_0+{\cal
L}(a_{22})d_{\nu}u_0+d^2u_0(\nu,\nu)\rho.$$ Straightforward
computation from the expressions for the $a_{ij}$ given in section 2
shows that:
$${\cal L}(a_{11})=k_0\rho,\quad {\cal
L}(a_{21})=-d_{\tau_0}\rho,\quad {\cal L}(a_{22})=0.$$ Substituting
in the expression for ${\cal L}(v^{-1})$, we find that the
coefficient of $d_{\tau_0}\rho$ vanishes, yielding:
$${\cal L}(\frac 1v)=-\frac 1{v_0^3}\{Du_0\cdot
DU+[ d_{\nu}u_0d^2u_0(\nu,\nu)+d_{\tau_0}u_0d^2u_0(\nu,\tau_0)
+(d_{\tau_0}u_0)^2k_0]\rho\}$$ Restoring the superscript $I$, we
write this (always at $\Gamma_0$) in the form:
$$\frac 1{v^I}=\frac 1{v_0^I}-\frac 1{(v_0^I)^3}Du_0^I\cdot
DU^I+\gamma_0^I\rho+E^{(1)I}(\rho,d_{\tau_0}\rho,DU^I)$$ where:
$$\gamma_0^I=-\frac
1{(v_0^I)^3}[d_{\nu}u_0^Id^2u_0^I(\nu,\nu)+d_{\tau_0}u_0^Id^2u_0^I(\nu,\tau_0)
+(d_{\tau_0}u_0^I)^2k_0]$$ and
$|E^{(1)I}|=O(\rho^2+(d_{\tau_0}\rho)^2+|DU^I|^2)$, with constants
depending only on the initial data.\vspace{.2cm}

The first angle condition $(BC^1)_{U^I}$ then becomes, with
$E^{(1)}=\sum_I'E^{(1)I}$:
$${\sum_I}'\frac {Du_0^I}{(v_0^I)^3}\cdot
DU^I-\rho{\sum_I}'\gamma_0^I=E^{(1)}(\rho,d_{\tau_0}\rho,DU).$$ Now
make the substitution $\rho\rightarrow \delta_0(U^2-U^1)$, valid on
$\Gamma_0$. We find:
$${\sum_I}'\frac
{Du_0^I}{(v_0^I)^3}\cdot
DU^I+\gamma_0^{(1)}(U^2-U^1)=E^{(1)}(\rho,d_{\tau_0}\rho,DU)_{|\rho=\delta_0(U^2-U^1)},$$
where:
$$\gamma_0^{(1)}:=-\delta_0
{\sum_I}'\gamma_0^I.$$\vspace{.2cm}

The second angle condition $(BC^2)_{U^I}$ has the form:
$${\sum_I}'\frac{d_{\nu}U^I}{v^I}+{\sum_I}'\frac{d_{\nu}u_0^I}{v^I}+\rho{\sum_I}'\frac{D^2u_0^I(\nu,\nu)}{v_0^I}=0,$$
or (using the linearization of $v^{-1}$ computed above):
$${\sum_I}'(\frac
1{v_0^I}d_{\nu}U^I-\frac{d_{\nu}u_0^I}{(v_0^I)^3}Du_0^I\cdot
DU^I)+({\sum_I}'(d_{\nu}u_0^I)\gamma_0^I
+\frac{d^2u_0^I(\nu,\nu)}{v_0^I})\rho=E^{(2)},$$ where
$|E^{(2)}(\rho,d_{\tau_0}\rho,DU)|=O(\rho^2+(d_{\tau_0}\rho)^2+|DU|^2)$.
Again making the substitution  $\rho\rightarrow \delta_0(U^2-U^1)$,
we find:
$${\sum_I}'(\frac
1{v_0^I}d_{\nu}U^I-\frac{d_{\nu}u_0^I}{(v_0^I)^3}Du_0^I\cdot DU^I)
+\gamma_0^{(2)}(U^2-U^1)=E^{(2)}(\rho,d_{\tau_0}\rho,DU)_{|\rho=\delta_0(U^2-U^1)},$$
where:
$$\gamma_0^{(2)}=\delta_0{\sum_I}'((d_{\nu}u_0^I)\gamma_0^I+\frac{d^2u_0^I(\nu,\nu)}{v_0^I}).$$
\vspace{.2cm}

In summary, the conjugation conditions may be written in linearized
form:
$$\left \{ \begin{array}{lr}
{\sum_I}'\frac{U^I}{v_0^I}=0&(LBC^0)_{U^I}\\
{\sum_I}'B_0^{(1)I}\cdot
DU^I+\gamma_0^{(1)}(U^2-U^1)=E^{(1)}(U,DU)&(LBC^1)_{U^I}\\
{\sum_I}'B_0^{(2)I}\cdot
DU^I+\gamma_0^{(2)}(U^2-U^1)=E^{(2)}(U,DU)&(LBC^2)_{U^I}
\end{array} \right .$$
(with zero initial conditions.)  Here we defined:
$$B_0^{(1)I}:=\frac{Du_0^I}{(v_0^I)^3},\quad B_0^{(2)I}:=\frac
1{v_0^I}\nu-\frac{d_{\nu}u_0^I}{(v_0^I)^3}Du_0^I.$$ Observe (for
future use) that, since
$1-(d_{\nu}u_0)^2/(v_0)^2=[1+(d_{\tau_0}u_0)^2]/(v_0)^2$, we have:
$${\sum_I}'B_0^{(2)I}\cdot
DU^I=[1+(d_{\tau_0}u_0)^2]{\sum_I}'\frac
1{(v_0^I)^3}d_{\nu}U^I-(d_{\tau_0}u_0){\sum_I}'\frac{d_{\nu}u_0^I}{(v_0^I)^3}d_{\tau_0}U^I.$$
(The notation uses the fact that $d_{\tau_0}u_0^I$ is independent of
$I$ at points of $\Gamma_0$).\vspace{.3cm}

\textbf{4.2 The fixed-point scheme.}\vspace{.2cm}

We write the quasilinear system $(PDE)_{U^I}/(BC)_{U^I}$ in a
slightly modified `linearized form', as follows:
$$(LPDE/LBC)_{U}\left \{ \begin{array}{l}
L_{h_0^I}[U^I]=f_U^I(x,t)\qquad\mbox{ on } D_0^I\times (0,T]\\
{\sum_I'}\frac{U^I}{v_0^I}=0,\quad {\mathbb
B}_0[DU]+\gamma_0[U]=E_U\mbox{ on }\Gamma_0,\quad
d_{\omega}U^3_{|\partial \Omega}=0 \end{array}\right. $$ where:
$$f^I_U(x,t):={\cal
F}^I_U+[(h_0^I)^{ab}-(h_U^I)^{ab}]\partial^2_{a,b}U^I-{\cal
A}^I_UL_{h^I_U}[{\cal E}[U^2-U^1]]$$ and we collapsed two linearized
angle conditions into one, introducing a slight change in notation.
$h_0^I=g_0^I$ is the induced metric at $t=0$.\vspace{.2cm}

This will be solved by a standard fixed-point argument in the space:
$${\mathbb X}_{R,T}:=\{U=(U^1,U^2,U^3); U^I\in
C^{2+\alpha,1+\alpha/2}(D_0^I\times [0,T]), U^I_{|t=0}\equiv 0\},$$
with  suitable choices of $R$ and $T$, to be described soon. (We
organize the fixed-point argument in the same way as
\cite{BaconneauLunardi}.)\vspace{.2cm}

Define an operator $\Phi:V\mapsto U$ by assigning to $V\in {\mathbb
X}_{R,T}$ the solution of the linear system:
$$(LPDE_0)\left \{\begin{array}{l}
L_{h_0^I}[U^I]=f_V^I(x,t)\qquad\mbox{ on } D_0^I\times (0,T]\\
{\sum_I'}\frac{U^I}{v_0^I}=0,\quad {\mathbb
B}_0[DU]+\gamma_0[U]=E_V\mbox{ on }\Gamma_0,\quad
d_{\omega}U^3_{|\partial \Omega}=0 \end{array}\right. $$

Consider the following assumptions for system $(LPDE_0)$:

(1) (Regularity of coefficients) $h_0^I\in C^{\alpha}(D_0^I)$,
${\mathbb B}_0,\gamma_0\in C^{1+\alpha}(\Gamma_0),(v_0^I)^{-1}\in
C^{2+\alpha}(\Gamma_0)$;

(2) (Complementarity) The conjugation operator
$({\sum_I}'\frac{U^I}{v_0^I},{\mathbb B}_0[DU])$ satisfies the
`complementarity conditions' on $\Gamma_0$, with respect to the
operators $L_{h^I_0}$.

Under these assumptions, it is a classical fact that this system
(with vanishing initial data) has a unique solution, provided
$f^I_V\in C^{\alpha,\alpha/2}(D_0^I\times [0,T]$ and $E_V\in
C^{1+\alpha,(1+\alpha)/2}(\Gamma_0\times [0,T])$ satisfy also:

(3) (Compatibility with $U^I_{|t=0}\equiv 0$:)
$${\sum_I}'
\frac {f^I_V}{v_0^I}_{|t=0}=0,\quad {E_V}_{|t=0}=0.$$

The solution satisfies the estimate:
$$||U^I||_{2+\alpha}\leq
M_0(\sum_I||f^I_V||_{\alpha}+||E_V||_{1+\alpha}),$$ with $M_0$
depending only on the H\"{o}lder norms of the coefficients in
(1).\vspace{.2cm}

The verification of condition (1) under the assumption $u_0^I\in
C^{3+\alpha}(D_0^I)$ is straightforward. Verifying (3) is also easy,
since when $U^I\equiv 0$ we have $E_U=0$ and
$f^I_U=tr_{h_0^I}D^2u_0^I$, so that the compatibility condition
amounts to:
$${\sum_I}'\frac 1{v_0^I}tr_{h_0^I}D^2u_0^I=0,$$
which is just the compatibility condition considered in section 1
(equivalent to ${\sum_I'}H_0^I=0$ on $\Gamma_0$.)

Complementarity (assumption (2)) also holds for the system
$(LPDE_0)$, but verifying this is more technical- it is done in the
next section. Assuming (2) for the moment, we conclude the map
$\Phi$ is well-defined.\vspace{.3cm}

\textbf{4.3 Contraction estimates in H\"{o}lder norms.}\vspace{.2cm}

To finish the argument, we must verify that (with suitable choices
of $R$ and $T$) $\Phi$ maps into ${\mathbb X}_{T,R}$ and is a
contraction in this space.\vspace{.2cm}

Let $V_1,V_2\in {\mathbb X}_{R,T}$. Then $W=\Phi(V_1)-\Phi(V_2)$ is
the solution of the linear problem:
$$(LPDE_W)\left \{\begin{array}{l}
L_{h_0^I}[W^I]=f_{V_1}^I-f_{V_2}^I\qquad\mbox{ on } D_0^I\times (0,T]\\
{\sum_I'}\frac{W^I}{v_0^I}=0,\quad {\mathbb
B}_0[DW]+\gamma_0[W]=E_{V_1}-E_{V_2}\mbox{ on }\Gamma_0,\quad
d_{\omega}W^3_{|\partial \Omega}=0 \end{array}\right. $$ (with zero
initial conditions.) Thus we have the estimate:
$$||W^I||_{2+\alpha}\leq
M_0(\sum_I||f_{V_1}^I-f_{V_2}^I||_{\alpha}+||E_{V_1}-E_{V_2}||_{1+\alpha}).$$
We write the difference $f_{V_1}^I-f_{V_2}^I$ as the sum of five
terms:
$$f_{V_1}^I-f_{V_2}^I={\cal F}^I_{V_1}-{\cal F}^I_{V_2}$$
$$+(h_{V_2}^{Iab}-h_{V_1}^{Iab})\partial^2_{ab}V_1^I+(h_0^{Iab}-h_{V_2}^{Iab})\partial^2_{ab}(V_1^I-V_2^I)$$
$$+({\cal A}^I_{V^1}-{\cal A}^I_{V_2})L_{h^I_{V_1}}[{\cal
E}(V_1^I)]$$
$$+{\cal
A}^I_{V^2}(h_{V_1}^{Iab}-h_{V_2}^{Iab})\partial^2_{ab}V_1^I$$
$$+{\cal A}_{V_2}^Ih_{V_2}^{Iab}\partial^2_{ab}(V_1^I-V_2^I)$$
$$:=f^{(1)I}+\ldots+f^{(5)I}.$$

We need bounds for each of the $f^{(i)I}$ in terms of
$||V_1-V_2||_{2+\alpha}$, $R$ and $T$. The details are lengthy but
standard, and it suffices to state the estimates with a brief
justification.($c_0$ denotes a constant depending only on the
initial data, which may change from one occurrence to the
next.)\vspace{.2cm}

 $$\quad ||{\cal F}^I_{V_1}-{\cal F}^I_{V_2}||_{\alpha}\leq
c_0T^{\alpha/2}||V_1-V_2||_{2+\alpha},$$ since ${\cal F}^I_V$
depends only on $V,DV,{\cal E}(V)$ and $D({\cal E}(V))$, and
vanishes at $t=0$. (The `nonlocality' of $\cal E$ does not make the
estimates harder, since this operator is linear in $V$ and bounded
in $C^{2+\alpha,1+\alpha/2}$ norm.)\vspace{.2cm}

$$||f^{(2)I}||_{\alpha}\leq
c_0T^{\alpha/2}||V_1-V_2||_{2+\alpha}R+c_0T^{\alpha/2}||V_1-V_2||_{2+\alpha},$$
since $h^I_{V_2}-h^I_{V_1}$and $h^I_0-h^I_{V_2}$ depend only on
derivatives of $V_1$, $V_2$ up to first order and vanish at
$t=0$.\vspace{.2cm}

$$({\cal A}^I_{V^1}-{\cal A}^I_{V_2})L_{h^I_{V_1}}[{\cal
E}(V_1^I)]\leq c_0T^{\alpha/2}||V_1-V_2||_{2+\alpha}R,$$ since
${\cal A}^I_{V_1}-{\cal A}^I_{V_2}$ depends on derivatives of
$V_1,V_2$ only up to first order and vanishes at $t=0$.
\vspace{.2cm}

$${\cal
A}^I_{V^2}(h_{V_1}^{Iab}-h_{V_2}^{Iab})\partial^2_{ab}V_1^I\leq
c_0T^{\alpha/2}||V_1-V_2||_{2+\alpha}R^2,$$ since ${\cal A}_{V_2}^I$
depends on derivatives of $V_2$ up to first order, and vanishes at
$t=0$.\vspace{.2cm}

Finally,
$${\cal A}_{V_2}^Ih_{V_2}^{Iab}\partial^2_{ab}(V_1^I-V_2^I)\leq
c_0T^{\alpha/2}||V_1-V_2||_{2+\alpha}R^2,$$ for the same
reason.\vspace{.2cm}

 A bit more involved (but also needed) is be the
estimate over $\Gamma_0$:
$$||E_{V_1}-E_{V_2}||_{1+\alpha}\leq
c_0T^{\alpha/2}||V_1-V_2||_{2+\alpha}R.$$ Although the $1+\alpha$
norm of $E_V$ involves $D^2V$, the estimate holds since $E_V$ is of
quadratic order in $V$ and $DV$.\vspace{.2cm}

We conclude that $\Phi$ is a 1/2 contraction, provided we choose $R$
and $T$ so that $16 c_0(1+R^2)T^{\alpha/2}\leq \frac 12$
($16=5\times 3+1$). Assuming this inequality holds, applying the
contraction estimate to the case $V_2=0$ yields, for any $V\in
{\mathbb X}_{R,T}$:
$$||\Phi(V)||_{2+\alpha}\leq ||\Phi(0)||_{2+\alpha}+\frac 12 R,$$
where $\Phi(0)$ is the solution (with zero initial data) of the
system $(LPDE)_0$ with data $f^I_V=tr_{h_0^I}D^2u_0$, $E_V=0$,
satisfying the estimate:
$$||\Phi(0)||_{2+\alpha}\leq M_0||tr_{h_0^I}D^2u_0^I||_{\alpha}.$$
\vspace{.2cm}

The choices of $R$ and $T$ are made as follows. As observed earlier,
if the $C^2$ geometry of the configuration at time $t$ is
sufficiently close to that at $t=0$, all the geometric constructions
used make sense; in particular the operator ${\cal E}$ and ${\cal
F}^I_U$ are well defined. Choose $r>0$ so that
$||U^I||_{C^2(D_0^I)}<r$ quantifies \ `$u^I$ sufficiently close to
$u_0^I$'. Then pick $T_0>0$ so that $U_{|t=0}\equiv 0$ and
$||U||_{2+\alpha}< T_0$ \emph{imply} $||U^I||_{C^2(D_0^I)}<r$. Now
choose $R>0$ sufficiently large to ensure:
$$ M_0||tr_{h_0^I}D^2u_0^I||_{\alpha}<\frac 12 R.$$
Finally, given $R$ we pick $T<T_0$ small enough that
$16c_0(1+R^2)T^{\alpha/2}\leq \frac 12$ holds.\vspace{.2cm}

With these choices we guarantee, on the one hand, that $\Phi$ is a
contraction; and, on the other, that:
$$||\Phi(V)||_{2+\alpha}\leq ||\Phi(0)||_{2+\alpha}+\frac 12 R <R$$
for any $V\in {\mathbb X}_{R,T}$. Thus $\Phi$ maps into $ {\mathbb
X}_{R,T}$, and has a unique fixed point $U$. This fixed point is the
unique solution of $(LPDE)_{U}/(LBC)_U$ with zero initial data. This
concludes the proof of local existence for the system
$(PDE)_{(u^I,\tilde{\rho})}$, and hence for the original system
$(PDE)_{w^I}$. (Except for the verification of complementarity for
$(LPDE)_0$, carried out in the next section).\vspace{.5cm}

\textbf{5. The complementarity condition.}\vspace{.2cm}

In this section we verify that the system of conjugation conditions
along $\Gamma_0$ satisfies the `complementarity conditions'
(Lopatinski-Shapiro) with respect to the $3\times 3$ linear
parabolic operator $L_{g_0^I}[U^I]$. (We use the conditions as
stated in \cite{EidelmanZhitarasu}.) Explicitly, the operator is:
$$L_{g_0^I}U^I=\partial_tU^I-g_0^{Iab}\partial^2_{ab}U^I\mbox{ on
}\Omega\times [0,T],$$
$$g_{0ab}^I=\delta_{ab}+\partial_au_0^I\partial_bu_0^I,\quad
g_0^{Iab}=\delta_{ab}-\frac{\partial_a u_0^I\partial_b
u_0^I}{(v_0^I)^2},$$
$$v_0^I=\sqrt{1+(d_{\tau_0}u_0^I)^2+(d_{\nu}u_0^I)^2}\mbox{ on
}\Gamma_0.$$ The conjugation operator has three components:
$$\left \{ \begin{array}{l}
B_0[Y]={\sum_I}'\frac{U^I}{v_0^I}\\
B_1[Y]={\sum_I}'\frac{n_0^I}{(v_0^I)^2}d_{\nu}U^I+\lambda_0{\sum_I}'\frac
1{(v_0)^2}d_{\tau_0}U^I,\\
B_2[Y]=(1+\lambda_0^2){\sum_I}'\frac
1{(v_0^I)^3}d_{\nu}U^I-\lambda_0{\sum_I}'\frac{n_0^I}{(v_0^I)^3}d_{\tau_0}U^I,
\end{array}\right .$$
where in this section we adopt the notation:
$$n_0^I:=d_{\nu}{u_0^I}_{|\Gamma_0},\quad
\lambda_0:=d_{\tau_0}{u_0^I}_{|\Gamma_0}.$$ Fix $x_0\in \Gamma_0$
and `straighten the boundary' via a diffeomorphism $\chi$ from a
neighborhood of $0\in {\mathbb R}^2$ (with coordinates $z=(z_1,z_2)$
to a neighborhood of $x_0$ in $\Omega$, and with the mapping
properties:
$$\chi:\{z_1=0\}\rightarrow \Gamma_0,\quad \chi(0)=x_0,\quad
d\chi(0):\partial_{z_1}\mapsto \nu(x_0),\quad
d\chi(0):\partial_{z_2}\mapsto \tau_0(x_0).$$

Let $Y^I(z,t)=U^I(\chi(z),t).$ From this point on the symbols
$\gamma^I$, $n_0^I$, $\lambda_0$, $v_0^I$ will denote the values of
the corresponding functions at the fixed point $x_0\in \Gamma_0$.

The transformed operator is:
$$L_0[Y^I]:=\partial_tY^I-\gamma^{Iij}\partial^2_{z_iz_j}Y^I.$$
 $\gamma^I=\chi(0)^*g_0^I(x_0)$ is the
pullback metric tensor, with components:
$$\gamma_{11}^I=1+(n_0^I)^2,\quad
\gamma_{12}^I=\lambda_0n_0^I,\quad \gamma_{22}^I=1+\lambda_0^2,$$
and inverse:
$$\gamma^{I11}=\frac{1+\lambda_0^2}{(v_0^I)^2},\quad
\gamma^{I12}=-\frac{\lambda_0n_0^I}{(v_0^I)^2},\quad
\gamma^{I22}=\frac{1+(n_0^I)^2}{(v_0^I)^2}.$$ In the new
coordinates, the components of the conjugation operator are:
$$\left \{ \begin{array}{l}
B_0[Y]={\sum_I}'\frac{Y^I}{v_0^I}\\
B_1[Y]={\sum_I}'\frac{n_0^I}{(v_0^I)^2}\partial_{z_1}Y^I+\lambda_0{\sum_I}'\frac
1{(v_0^I)^2}\partial_{z_2}Y^I,\\
B_2[Y]=(1+\lambda_0^2){\sum_I}'\frac
1{(v_0^I)^3}\partial_{z_1}Y^I-\lambda_0{\sum_I}'\frac{n_0^I}{(v_0^I)^3}\partial_{z_2}Y^I,
\end{array}\right .$$\vspace{.2cm}

The next step is to consider the Fourier-Laplace transform
$\hat{Y}^I(\tau,\xi,p)$ of $Y^I(z_1,z_2,t)$: Fourier transform in
$z_2$ corresponds to the variable $\xi\in \mathbb{R}$, while Laplace
transform in $t$ corresponds to $p\in \mathbb{C}$. Furthermore, we
adjust the signs so that $\{\tau>0\}$ corresponds to $D_0^1=D_0^2$
(that is, $\tau=z_1$ for $\hat{Y}^1$ and $\hat{Y}^2$) and also to
$D_0^3$ (that is, $\tau=-z_1$ for $\hat{Y}^3$). (This is how we deal
with the fact that we have a conjugation problem, rather than a
boundary-value problem.) This introduces a sign $\sigma^I$ in the
transformed operator ($\sigma^1=\sigma^2=1,\sigma^3=-1$), which is
an ordinary differential operator in the variable $\tau$, for each
fixed $(p,\xi)$:
$$\gamma^{I11}\frac
{d^2\hat{Y}^I}{d\tau^2}+2\sigma^Ii\xi\gamma^{I12}\frac{d\hat{Y}^I}{d\tau}-\xi^2\gamma^{I22}\hat{Y}^I
-p\hat{Y}^I.$$ We are interested in solutions of the corresponding
ODE which decay as $\tau \rightarrow +\infty$. If we assume the
exponential form: $$\hat{Y}^I_{\xi,p}(\tau)=e^{i\tau
\rho_I}\hat{Y}^I_{\xi,p}(0),$$ this corresponds to
$\mathbb{I}m(\rho_I)>0$. The $\rho_I\in \mathbb{C}$ depend on $\xi$
and $p$, and are roots of the indicial equation:
$$\gamma^{I11}\rho_I^2+2\sigma^I\gamma^{I12}\xi\rho_I+\gamma^{I22}\xi^2+p=0.$$
Using the expressions given earlier for the $\gamma^{Iij}$, we find
for the roots:
$$\rho_I=\frac
1{\gamma^{I11}}[-\sigma^I\gamma^{I12}\xi+i\sqrt{\xi^2\det[\gamma^I]^{-1}+p\gamma^{I11}}]$$
$$=\frac{v_0^I}{1+\lambda_0^2}[\sigma^I\lambda_0\frac{n_0^I}{v_0^I}\xi+i\sqrt{\Delta}],$$
where $\Delta=\xi^2+p(1+\lambda_0^2)$ and we take the branch
${\mathbb R}e(\sqrt{\Delta})>0$. The boundary operators for the ODE
are obtained from the correspondence:
$$\partial_{z_1}Y^I_{|z_1=0}\longrightarrow
\sigma^I\frac{d\hat{Y}^I}{d\tau}_{|\tau=0}=\sigma^Ii\rho_I\hat{Y}^I(0),\quad
\partial_{z_2}Y^I_{|z_1=0}\longrightarrow i\xi\hat{Y}^I(0).$$
We find:
$$\begin{array}{l}
B_0[Y]\longrightarrow {\sum_I}'\frac{\hat{Y}^I(0)}{v_0^I}\\
B_1[Y]\longrightarrow
{\sum_I}'[\frac{n_0^I}{(v_0^I)^3}\sigma^I(i\rho_I)+\frac{\lambda_0}{(v_0^I)^3}(i\xi)]\hat{Y}^I(0)\\
B_2[Y]\longrightarrow
{\sum_I}'[\frac{1+\lambda_0^2}{(v_0^I)^3}\sigma^I(i\rho_I)-\frac{\lambda_0n_0^I}{(v_0^I)^3}(i\xi)]\hat{Y}^I(0)
\end{array}$$\vspace{.2cm}

The complementarity condition is the statement that, for each
$(p,\xi)$ in the set:
$${\cal A}:=\{(p,\xi)\in {\mathbb C}\times
\mathbb{R};|p|+|\xi|>0,(1+\lambda_0^2)\mathbb{R}e(p)>-\xi^2\}$$ the
map $\hat{Y}^I(0)\mapsto B_I[Y]$ defines an isomorphism of
$\mathbb{C}^3$. In matrix form (after canceling a common factor
$(v_0^I)^{-1}$ in each column) we are interested in:
$$\mathbb{T}_{(p,\xi)}:=\left [ \begin {array}{ccc}
1&1&-1\\\frac{n_0^1}{(v_0^1)^2}\rho_1+\frac{\lambda_0}{(v_0^1)^2}\xi&
\frac{n_0^2}{(v_0^2)^2}\rho_2+\frac{\lambda_0}{(v_0^2)^2}\xi&
\frac{n_0^3}{(v_0^3)^2}\rho_3-\frac{\lambda_0}{(v_0^3)^2}\xi\\
\frac{1+\lambda_0^2}{(v_0^1)^2}\rho_1-\frac{\lambda_0n_0^1}{(v_0^1)^2}\xi&
\frac{1+\lambda_0^2}{(v_0^2)^2}\rho_2-\frac{\lambda_0n_0^2}{(v_0^2)^2}\xi&
\frac{1+\lambda_0^2}{(v_0^3)^2}\rho_3+\frac{\lambda_0n_0^3}{(v_0^3)^2}\xi
\end{array} \right ]$$\vspace{.2cm}

\textbf{Lemma 5.1.} For each $(p,\xi)\in {\cal A}$,
$\det(\mathbb{T}_{(p,\xi)})\neq 0$.\vspace{.3cm}

\emph{Proof.} We begin by introducing some notation. Consider the
vectors $t_I\in \mathbb{R}^2,r_I\in \mathbb{C}^2$:
$$t_I=\frac 1{v_0^I}[n_0^I,\sqrt{1+\lambda_0^2}],\quad r_I=\frac
1{v_0^I}[\sqrt{1+\lambda_0^2}\rho_I,\lambda_0\sigma^I\xi].$$ Then:
$$\frac{n_0^1}{v_0^1}\frac{\rho_1}{v_0^1}+\frac{\lambda_0}{v_0^1}\frac{\xi}{v_0^1}=\frac
1{\sqrt{1+\lambda_0^2}}t_1\cdot r_1,\quad
\frac{1+\lambda_0^2}{v_0^1}\frac{\rho_1}{v_0^1}-\frac{n_0^1}{v_0^1}\frac{\lambda_0\xi}{v_0^1}=t_1\cdot
r_1^{\perp},$$ and analogously for the other two columns. (Here we
use the `real dot product', i.e. $t_I\cdot r_I=t_I\cdot
\mathbb{R}e(r_I)+it_I\cdot \mathbb{I}m(r_I)$).\vspace{.2cm}

With this notation, we find (after cancelation of common factors and
using $t_I\cdot r_I^{\perp}=-t_I^{\perp}\cdot r_I$) we must show
that the matrix:
$$\left [ \begin{array}{ccc}
1&1&-1\\
t_1\cdot r_1&t_2\cdot r_2&t_3\cdot r_3\\
t_1^{\perp}\cdot r_1&t_2^{\perp}\cdot r_2&t_3^{\perp}\cdot r_3
\end{array}\right ]$$
has non-zero determinant \textbf{det}. Expanding, we have:
$$\mbox{\textbf{det}}=(t_1^{\perp}\cdot r_1)(t_2\cdot r_2)-(t_2^{\perp}\cdot
r_2)(t_1\cdot r_1)-(t_1\cdot r_1)(t_3^{\perp}\cdot
r_3)+(t_1^{\perp}\cdot r_1)(t_3\cdot r_3)+(t_2\cdot
r_2)(t_3^{\perp}\cdot r_3)-(t_2^{\perp}\cdot r_2)(t_3\cdot
r_3).$$\vspace{.2cm}

The argument has two parts. In the first part, we use the symmetries
in the conjugation conditions, encoded in the properties of $t_I$:
$$t_1+t_2=t_3,\quad |t_I|=1.$$
From the geometric version of lemma 1.2, we can express $t_1$ and
$t_2$ in terms of $t_3$:
$$t_1=R_{\pi/3}t_3=\frac 12 t_3+\frac{\sqrt{3}}2t_3^{\perp},\quad
t_2=R_{5\pi/3}t_3=\frac 12 t_3-\frac{\sqrt{3}}2t_3^{\perp}.$$
Substituting these values in the expression for \textbf{det}, we
obtain:
$$\mbox{\textbf{det}}=\frac{\sqrt{3}}2(r_1+r_2)\cdot[(t_3\cdot
r_3)t_3+(t_3^{\perp}\cdot r_3)t_3^{\perp}]+\frac
12(r_1-r_2)\cdot[(t_3^{\perp}\cdot r_3)t_3-(t_3\cdot
r_3)t_3^{\perp}]$$
$$+r_2\cdot \{\frac{\sqrt{3}}2[(r_1\cdot t_3)t_3+(r_1\cdot
t_3^{\perp})t_3^{\perp}]+\frac 12[(t_3^{\perp}\cdot
r_1)t_3-(t_3\cdot r_1)t_3^{\perp}]\},$$ which simplifies to give an
expression in terms of the $r_I$ only:
$$\mbox{\textbf{det}}=\frac{\sqrt{3}}2(r_1\cdot r_3+r_2\cdot
r_3+r_1\cdot r_2)+\frac 12(r_1\cdot r_2^{\perp}+r_2\cdot
r_3^{\perp}+r_3\cdot r_1^{\perp}).$$\vspace{.2cm}

In the second part of the argument, we use symmetries in the $r_I$.
Recalling the expressions for the roots $\rho_I$, we have (with
$a_I:=\mathbb{R}e(r_I),b_I:=\mathbb{I}m(r_I)$, both in
$\mathbb{R}^2$):
$$r_I=a_I+ib_I,\quad a_I=\lambda_0\xi[\frac
1{\sqrt{1+\lambda_0^2}}\frac{\sigma^In_0^I}{v_0^I},\frac{\sigma^I}{v_0^I}],
\quad b_I=\frac{\sqrt{\Delta}}{\sqrt{1+\lambda_0^2}}[1,0].$$ We see
that $b_I:=b$ is independent of $I$, while the $a_I$ satisfy:
$$\sum_Ia_I=0,\quad
|a_I|^2=\frac{\lambda_0^2\xi^2}{1+\lambda_0^2}:=|a|^2$$ (independent
of $I$). Consider first the imaginary part of \textbf{det}. We have:
$$\mathbb{I}m(r_I\cdot
r_J)=\mathbb{R}e(r_I)\mathbb{I}m(r_J)+\mathbb{R}e(r_J)\mathbb{I}m(r_I)=(a_I+a_J)\cdot
b,$$ and from this it follows easily that $\mathbb{I}m(r_1\cdot
r_3+r_2\cdot r_3+r_1\cdot r_2)=0$. Similarly, since
$r_I^{\perp}=a_I^{\perp}+ib^{\perp}$, we have:
$$\mathbb{I}m(r_I\cdot r_J^{\perp})=a_I\cdot
b^{\perp}+a_J^{\perp}\cdot b,$$ and again we have:
$\mathbb{I}m(r_1\cdot r_2^{\perp}+r_2\cdot r_3^{\perp}+r_3\cdot
r_1^{\perp})=0$. Thus $\mathbb{I}m$(\textbf{det})$=0$.\vspace{.2cm}

Turning to the real part, since $\mathbb{I}m(r_I)\cdot
\mathbb{I}m(r_I)^{\perp}=b\cdot b^{\perp}=0$:
$$\mathbb{R}e(r_I\cdot r_J)=a_I\cdot a_J-|b|^2,\quad
\mathbb{R}e(r_I\cdot r_J^{\perp})=a_I\cdot a_J^{\perp}-b\cdot
b^{\perp}=a_I\cdot a_J^{\perp}.$$

Using the version of lemma 1.2 for the case $\sum_Ia_I=0,|a_I|=|a|$:
$$a_1=R_{4\pi/3}a_3,\quad a_2=R_{2\pi/3}a_3$$
we see that:
$$a_1\cdot a_3=a_2\cdot a_3=a_1\cdot a_2=-\frac 12 |a|^2,$$
$$a_1\cdot a_2^{\perp}=a_2\cdot a_3^{\perp}=a_3\cdot
a_1^{\perp}=\frac{\sqrt{3}}2|a|^2.$$ We conclude:
$$\mbox{\textbf{det}}=\frac{\sqrt{3}}2(-\frac 32|a|^2-3|b|^2)+\frac
12 \frac
{3\sqrt{3}}2|a|^2=-(3\sqrt{3}/2)|b|^2=-(3\sqrt{3}/2)\frac{|\Delta|}{1+\lambda_0^2}.$$
This is non-zero for any $(p,\xi)\in {\cal A}$ (since
$\mathbb{R}e(\Delta)>0$ in $\cal A$), as we had to
show.\vspace{.2cm}

\textbf{Remark.} Note that this calculation depends on the fact that
the boundary conditions (that is, the $t_I$) and the roots of the
characteristic equation (that is, the $r_I$, and ultimately the
coefficients of the operator, which come from the induced metric at
the junction) satisfy the same type of symmetry. That is,
complementarity of this particular set of conjugation conditions
seems to be linked to the fact that we are dealing with mean
curvature motion; this suggests that local well-posedness may fail
for a more general parabolic system (with the same conjugation
conditions.)


\newpage

\begin{thebibliography}{99}

\bibitem[BaconneauLunardi]{BaconneauLunardi} Baconneau, O.; Lunardi,
A. \emph{ Smooth solutions to a class of free boundary parabolic
problems}.  Trans. Amer. Math. Soc.  \textbf{356}  (2004),  no. 3,
987--1005.

\bibitem[BronsardReitich]{BronsardReitich}Bronsard, Lia; Reitich, Fernando
 On three-phase boundary motion and the singular limit
 of a vector-valued Ginzburg-Landau equation.  Arch. Rational Mech. Anal.\textbf{124}  (1993),
   no. 4, 355--379.

\bibitem[Dierkes et al.]{Dierkes et al.}U.Dierkes,
S.Hildebrandt A.K\"{u}ster, O.Wohlrab, \emph{Minimal Surfaces 1.}
Grunlehren der Mathematischen Wissenschaften v.295, Springer-Verlag
1992. ISBN 3-540-53169-6

\bibitem[EidelmanZhitarasu]{EidelmanZhitarasu}
Eidelman, S. D.; Zhitarashu, N. V. \emph{Parabolic boundary value
problems.}  Operator Theory: Advances and Applications, 101.
Birkhäuser Verlag, Basel, 1998. xii+298 pp. ISBN: 3-7643-2972-6

\bibitem[Freire]{Freire} Freire,A. \emph{Mean Curvature Motion of Graphs with
Constant Contact Angle and Moving Boundaries} (preprint, May
2008-arXiv:0805.4592)

\bibitem[LensSeminar]{LensSeminar}
 \emph{Evolution of convex lens-shaped networks under curve shortening flow}
 O. Schn\"{u}rer, A. Azouani, M. Georgi, J. Hell, N. Jangle, A. Koeller,
 T. Marxen, S. Ritthaler, M. S\'{a}ez, F. Schulze, B. Smith,
 (Lens Seminar, FU Berlin 2007), arXiv:0711.1108


\bibitem[Mantegazza et al.]{Mantegazza et al.} Mantegazza,C.,
Novaga,M., Tortorelli,V., \emph{Motion by curvature of planar
networks.} Ann. Sc. Norm. Super. Pisa Cl. Sci. \textbf{5}(3) (2004)
no. 2, 235--324.


\end{thebibliography}
\end{document}